\documentclass[11pt]{amsart}
\usepackage{amsfonts,amscd}
\usepackage{epsfig}
\def\N{{\mathbb N}}
\def\Z{{\mathbb Z}}

\def\R{{\mathbb R}}

\def\P{{\mathbb P}}
\def\gp{\gamma_{p}}
\def\gpd{\gamma_{p,q}}
\def\gm{\gamma_{M}}
\newtheorem{thr}{Theorem}[section]

\newtheorem{lm}[thr]{Lemma}
\newtheorem{pr}[thr]{Proposition}
\newtheorem{co}[thr]{Corollary}

\newtheorem{prob}[thr]{Problem}
\newtheorem{thr*}{Theorem}
\newtheorem{pr*}{Proposition}
\newtheorem{co*}{Corollary}
\begin{document} 
\title[Lagrangian Surfaces in a Fixed Homology Class] 
{Lagrangian Surfaces in a Fixed Homology Class: \\
Existence of Knotted Lagrangian Tori}
\author{Stefano Vidussi}
\address{Department of Mathematics, University of California Riverside, Riverside, CA 92521}
\email{svidussi@math.ucr.edu}
\maketitle
\baselineskip 18pt
\noindent {\bf Abstract.} In this paper we show that there exist
simply connected symplectic manifolds which contain infinitely many
knotted lagrangian tori, i.e., nonisotopic lagrangian tori that are image of homotopic
embeddings. Moreover, the homology class they represent can be
assumed to be nontrivial and primitive. This answers a question of
Eliashberg and Polterovich.
\section{Introduction and statement of the results}
Let $M$ be a smooth, closed $4$-manifold endowed with a symplectic form $\omega$
and let $\Sigma$ be a smooth, closed oriented
surface. Consider two symplectic (resp. lagrangian) embeddings
$\phi_{1}$ and $\phi_{2}$ of $\Sigma$ in $M$.
Assume furthermore that the $\phi_{i}$'s are smoothly homotopic.
Among these maps we can define two
notions of equivalence, the second implying the first:
\vspace*{6pt}

\noindent (E1) The $\phi_{i}$ are {\it smoothly} isotopic, i.e., there exists a
smooth homotopy composed of embeddings;

\noindent (E2) The $\phi_{i}$ are {\it symplectically}
isotopic, i.e., there exists a smooth isotopy composed of symplectic
(resp. lagrangian) embeddings.
\vspace*{6pt} \\  (In what follows, the terms {\it homotopy} and {\it isotopy}
will always refer to smooth ones.)

A necessary condition for two embeddings to be equivalent under the equivalence
relations above is that their images, 
the embedded surfaces $\Sigma_{i} := \phi_{i}(\Sigma)$,
satisfy the corresponding equivalence relation (that will be referred to with 
the same notation): in the case of (E1),
the surfaces $\Sigma_{i}$ must be isotopic submanifolds of $M$, while
in case of (E2), the surfaces $\Sigma_{i}$ must be isotopic through symplectic (resp. lagrangian) submanifolds.
(Note that if the genus of $\Sigma$ is greater than $0$, these conditions could be not sufficient,
as the embeddings $\phi$ and $\phi \cdot \gamma$, where $\gamma$ a selfdiffeomorphism of $\Sigma$
that is not isotopic to the identity, have the same image but could be nonisotopic.)

In what follows we will concentrate on the isotopy problem for the surfaces $\Sigma_{i}$;
observe that, by a standard argument, two embeddings of a surface in a simply connected
$4$-manifold $M$ are homotopic if and only if their images represent the same homology class.

A priori, the first equivalence relation belongs to
differential topology, while the second one belongs to symplectic
topology. However, as we assume that the embeddings are symplectic or
lagrangian, in understanding (E1) we have to take
into account the constraint related to the rigidity induced by this extra
condition. In particular, this could prevent us from the possibility of
realizing an equivalence class of embeddings with a symplectic or lagrangian representative. As we will see in the following, this rigidity can affect the
existence of different classes of embeddings modulo the equivalence (E1).

The classification of embedded surfaces, modulo one of the equivalence relations
above, can be defined as the {\it symplectic} (resp. {\it lagrangian}) {\it knot
problem}. In particular, homologous but nonisotopic surfaces determine
different {\it knotted} (in the sense of differential topology) symplectic or
lagrangian surfaces, while isotopic surfaces that are not isotopic in the
symplectic sense determine different {\it symplectically knotted} symplectic
or lagrangian surfaces. (Note that, at least in general, there is no
``unknot'' i.e., a preferred representative of an homology class.)

In the recent past, a series of papers has addressed the question of
determining whether two embedded symplectic surfaces representing the same
homology class in a simply connected symplectic $4$-manifold $M$ must be
isotopic. One motivation for the isotopy problem for symplectic manifolds
comes from the analogous question for the case of complex curves on K\"ahler surfaces, where it is known, by classical results,
that complex representatives of the same homology class are in fact
isotopic. Working in this direction, Siebert and Tian have proven
that a symplectic surface $\Sigma_{d}$, representing the class
$d[H]$ in $H^{2}(\P^{2},\Z)$
(with standard notation) is symplectically isotopic to an algebraic curve of
degree
$d$, at least for $d \leq 17$. (Note that this result is stronger that the
one holding in the K\"ahler case, as now we are requiring only that
$\Sigma_{d}$ is symplectic.) A similar result holds for surfaces in $S^{2}
\times S^{2}$. These results are presented in \cite{Ti}.

In contrast with that, Fintushel and Stern have presented
a method to build nonisotopic symplectic tori representing
multiples of the class of a symplectic $c$-embedded torus in
a symplectic $4$-manifold (see \cite{FS2} for precise definitions and results).
Their construction, that applies to a large class of symplectic
manifolds, shows in particular the
existence of infinitely many nonisotopic simplectic tori representing
the homology class $2m[F]$ ($m \geq 2$) for
any elliptic surface $E(n)$ of fiber $F$. The latter result has been extended
by the author in \cite{V} to cover the case of all multiples $q[F]$
(at least for $n \geq 3$); in particular, when $q=1$, we see that there exist
infinitely many symplectic surfaces homologous, but nonisotopic, to
a complex connected curve ($F$ itself). Etg\"u and Park have presented
in \cite{EP} further constructions of nonisotopic tori
(in classes with divisibility) that complement the previous ones.
Examples of classes with positive
self-intersection, or higher genus, are still eluding us (at least for simply
connected manifolds: otherwise, see \cite{S}). Some of the non-isotopy results
above have been analyzed by Auroux, Donaldson and Katzarkov in \cite{ADK},
where they show that these examples (and new ones they built for the case of
singular symplectic surfaces in $\P^{2}$) could be interpreted as a kind of
braiding of parallel copies of a complex curve. The openness of the symplectic
condition allows us to keep the submanifolds resulting from this braiding
symplectic.

Considered all together, these results imply that, in suitable
manifolds, infinitely many knot types can be realized by symplectic surfaces.
On the other hand, the author is not aware of any example of isotopic
symplectic surfaces that are symplectically knotted.

The symplectic knot problem, however, has been preceded historically by
its lagrangian counterpart that, for reasons detailed below, is somewhat more
subtle. This question arose, for the case of
lagrangian $\R^{n}$'s in $\R^{2n}$ linear outside a ball, in the
``first paper'' on symplectic topology by Arnol'd (see \cite{A}).
In a more general set up, the problem has been summarized by Eliashberg and
Polterovich in \cite{ElP3}.

The analogy between symplectic and lagrangian knot problems is obvious.
But there are reasons that make the latter question more subtle, at least
in relation to the equivalence (E1).
The first (and probably less relevant) is that we can perturb the
symplectic form $\omega$ on $M$ to a symplectic form $\omega_{\epsilon}$
in such a way that an essential surface $\Sigma$, lagrangian with respect to
$\omega$, is symplectic with respect to $\omega_{\epsilon}$.
This result, whose proof appears in Gompf (\cite{G}), holds true also for
pairs of disjoint surfaces, so that the existence of essential knotted
lagrangian surfaces entails the existence of knotted symplectic surfaces.
But the main reason of interest stems from the fact that the lagrangian
condition is a closed one, in that respect similar to the condition of
being complex. In particular the rigidity of this condition imposes
constraints to the possibility of braiding copies of
a lagrangian surface in the spirit of \cite{ADK}, and makes a result of existence
of lagrangian knots appear more problematic.

In fact, the few results known so far point towards absence of knotted
lagrangian surfaces.
In \cite{ElP1} Eliashberg and Polterovich show that when $\Sigma = S^{2}$ or
$T^{2}$ a lagrangian surface in $T^{*}\Sigma$ homologous to the zero section
is in fact isotopic to the zero section (this result is quite exhaustive as,
by \cite{LS}, every homologically nontrivial lagrangian submanifold of
$T^{*}\Sigma$ must be homologous to the zero section). With similar spirit,
the authors prove that, at a local level (see \cite{ElP2} for precise
definitions and statements) lagrangian submanifolds are unknotted (in the
sense of (E2)). Again, in \cite{L} Luttinger proved that infinitely many knot types
of tori in $\R^{4}$ can not be realized with lagrangian embeddings.
In light of these results, Question 1.3A of \cite{ElP2}
asks whether there exist homotopic, but not symplectically isotopic, embeddings
of a lagrangian surface. Clearly, examples of this kind are provided by 
homologous 
lagrangian surfaces that fail to satisfy (E1) or (E2).
Seidel has answered in the positive to this question
(see \cite{Se1} and \cite{Se2}) constructing an infinite number of lagrangian
spheres, in a suitable $4$-manifold,
that are isotopic but symplectically knotted, i.e., equivalent under (E1) but
not under (E2).

The goal of this paper is to complete the answer to
Question 1.3A of \cite{ElP2}
constructing an infinite number of knotted lagrangian tori, i.e., homologous lagrangian
tori that are not equivalent under (E1).
We will prove the following results:
\begin{thr} \label{infinite} Let $E(2)_{K}$ be the symplectic $4$-manifold
(homotopy equivalent to $E(2)$) obtained by knot surgery on the left-handed trefoil $K$;
then there exists a nontrivial primitive homology class $[R]$ such that
any multiple $q[R]$, $q \geq 1$, can be represented
by infinitely many mutually nonisotopic lagrangian tori. \end{thr}

Theorem \ref{infinite} asserts therefore that infinitely many knot types can
be realized by lagrangian embeddings.

It will be apparent from the proof that we are able in fact to prove something
stronger, namely that if we denote by $R_{i}$ the lagrangian tori of Theorem
\ref{infinite}, there is no diffeomorphism of pairs $(E(2)_{K},R_{i}) \rightarrow
(E(2)_{K},R_{j})$, even not connected to the identity, unless $i=j$.

This result can be extended without any effort to cover other classes of
symplectic knot surgery manifolds, using for example the figure-eight knot,
or any non-prime fibered knot containing $K$ as a
summand, and many others. In fact,
we expect the result to hold for all symplectic manifolds
obtained by knot surgery on a fibered nontrivial knot, although the
proof of this would require some modification in the proof of Theorem
\ref{infinite}. We will address this problem in the future.

Using the aforementioned result of Gompf, we obtain the following Corollary:
\begin{co} \label{small} The manifold $E(2)_{K}$ has a primitive homology
class which can be represented by mutually nonisotopic symplectic tori.
\end{co}
We point out that phenomena of this kind for $E(2)_{K}$ are discussed
in \cite{FS2} and \cite{EP}, but for homology classes with divisibility.
Corollary
\ref{small} gives therefore the simplest (in the sense of geography)
known symplectic manifold having symplectic nonisotopic tori
in the same primitive homology class. Obviously, the same result obtained
in Corollary \ref{small} holds for all multiples of $[R]$.

We will briefly overview the ideas underlying the proof of Theorem
\ref{infinite}. We will present the manifold $E(2)_{K}$ as
result of link surgery over the link $L = K \cup M$ given by the knot $K$
and its meridian $M$. The link exterior $S^{3} \setminus \nu L$
fibers over $S^{1}$ with fiber $\Sigma_{L}$ given by the fiber
$\Sigma_{K}$ of $K$ with a disk removed. We can obtain symplectic tori
by looking at curves in $S^{3} \setminus \nu L$ transverse to $\Sigma_{L}$.
This is the approach of all available constructions of symplectic knots.
Here, we will proceed instead in the opposite way:
Let $\gamma$ be a curve in $\Sigma_{L}$. Denote by $N_{K}$ the
fibered $3$-manifold obtained by surgery on $K \cup M$ with coefficients
$0$ and $\infty$ respectively. Then the torus $S^{1} \times \gamma$ is
a lagrangian, framed torus in the manifold $S^{1} \times N_{K}$ (endowed of
a standard symplectic structure). By symplectic fiber summing two copies of
$E(1)$ to $S^{1} \times N_{K}$, this torus defines a framed lagrangian torus
in $E(2)_{K}$. The problem at this point is reduced
to find infinitely many curves $\gamma \in \Sigma_{L}$ homologous (but
nonisotopic) in
$S^{3} \setminus \nu L$, and then try to distinguish the isotopy class of the
resulting lagrangian tori.
The latter result will be obtained with the technique we introduced in \cite{V},
namely by studying the Seiberg-Witten polynomial
of the (symplectic) $4$-manifolds given by fiber summing $E(1)$ along
the lagrangian tori.
As the sum with $E(1)$ does not depend on the choice of the gluing
map, or the framing, the smooth structure of the resulting manifolds is
determined by the smooth isotopy class of the tori, and SW theory allows us
to distinguish the manifolds to the degree required in Theorem
\ref{infinite}.

We finish by observing that the results above can be generalized without effort
to the symplectic manifold $E(n)_{K}$ for any $n \geq 2$. The case of $n=1$,
where the situation is somewhat different, will be discussed in future
research.

\noindent {\it Added in proof:} R. Fintushel and R. Stern have in fact extended Theorem \ref{infinite}
to all $E(n)_K$, where $K$ is any nontrivial fibered knot. Also, they show that for $n = 1$ the nullhomologous lagrangian 
tori resulting from our (and their) construction are not isotopic. See R. Fintushel, R. Stern, \textit{Invariants for Lagrangian tori}, Geom. Topol.,  8  (2004), 947--968.

\noindent {\bf Acknowledgements:} It's about time to pay my debt of gratitude
to Ron Fintushel and Ron Stern for their invaluable support in this
and other papers of mine. I would like to thank also Slaven Jabuka for a 
conversation that led me to investigate this problem, and Paul Seidel for
several discussions.

\section{Preliminaries} \label{prel}
In this Section we will recall some standard definitions and results that can be
found, for example, in \cite{GS}; we will moreover specify the different notions of
knot.

Let $M$ be a smooth, closed, simply connected $4$-manifold and $\Sigma$ a smooth, closed, oriented surface. Given a homotopy class of maps $h \in [\Sigma,M]$, we will say that
two embeddings $\phi_{1},\phi_{2} \in h$ are isotopic
(equivalence relation (E1) ) if there is a
homotopy $\phi_{t} : \Sigma \rightarrow M$ through embeddings.
The images $\Sigma_{i} := \phi_{i}(\Sigma)$ of isotopic embeddings are isotopic submanifolds of $M$, i.e., there is an isotopy of the inclusion map of the first one
that has as image of the terminal map the second.
By the Isotopy Extension Theorem, isotopic submanifolds are
ambient isotopic, i.e., there exists a self-diffeomorphism $f$ of $M$ connected
to the identity such that $f(\Sigma_{1}) = \Sigma_{2}$.
Because of that, when we glue a manifold to $M$ along a submanifold $N$,
the result of the surgery depends, up to diffeomorphism, only on the
isotopy class of $N$, together with the choice of a gluing map on $N$ and 
of a lifting of this map to $\partial \nu N$.
 
We will say that $M$ contains a knotted surface $\Sigma$
if there exists an homotopy class of maps in
$[\Sigma,M]$ whose images represent (infinitely many) isotopy classes of embedded surfaces.
By standard arguments, this corresponds to have nonisotopic
representatives of the same homology
class of $M$. Each isotopy class is called a knot, or a knot type.

When $M$ admits a symplectic structure $\omega$, an embedding $\phi : \Sigma
\rightarrow M$ is lagrangian if $\phi^{*} \omega$ in trivial on $\Sigma$
or, which is the same, it restricts trivially to the image $\phi(\Sigma)$.
Coherently with the previous definition, we say that $M$ contains a knotted lagrangian 
surface if, for an homotopy class in $[\Sigma,M]$, (infinitely many) knot types are realized by embedded lagrangian surfaces. 

We remark that we could define the equivalence of embedded surfaces in terms of the coarser definition of pair diffeomorphism instead than isotopy (for classical knots in $S^{3}$ the definitions coincide). As previously mentioned, our results hold also under this more restrictive condition. 

As observed before, we could further ask whether
the isotopy between lagrangian surfaces can be realized through lagrangian
surfaces (equivalence relation E(2)). We say that $M$ contains a symplectically knotted surface $\Sigma$ if there are (infinitely many) lagrangian embeddings of $\Sigma$ whose images are isotopic, but not
symplectically isotopic. See \cite{Se1} and \cite{Se2} for the results on this.

\section{Construction of the links} \label{linkco}
In this section we will discuss  the construction of  a family of links, that
will be useful in the proof of Theorem \ref{infinite}.

Let $K$ be a nontrivial fibered knot, and denote by $L$ the
$2$-component link given by $K$ itself and a meridian $M$. Consider
the link exterior $S^{3} \setminus \nu L$. 
We have the following simple Proposition, that follows easily by observing 
that, up to isotopy, we can assume that $M$ is transverse to the fiber
$\Sigma_{K}$ of $S^{3} \setminus \nu K$.
\begin{pr} The $3$-manifold $S^{3} \setminus \nu L$ admits a fibration 
over $S^{1}$, having class $(1,0) \in H^{1}(S^{3} \setminus \nu L,\Z) = 
\Z^{2}$, with fiber $\Sigma_{L}$ given by the spanning surface
$\Sigma_{K}$ of $K$ with a disk
removed. \end{pr} In the Proposition above we have implicitly assumed as 
cohomology basis the basis of $H^{1}(S^{3} \setminus \nu L,\Z)$ dual to
$[\mu(K)],[\mu(M)]$. We want to remark that $\Sigma_{L}$ {\it is not} the 
spanning surface of $L$. 

\begin{figure}[h] 
\centerline{\psfig{figure=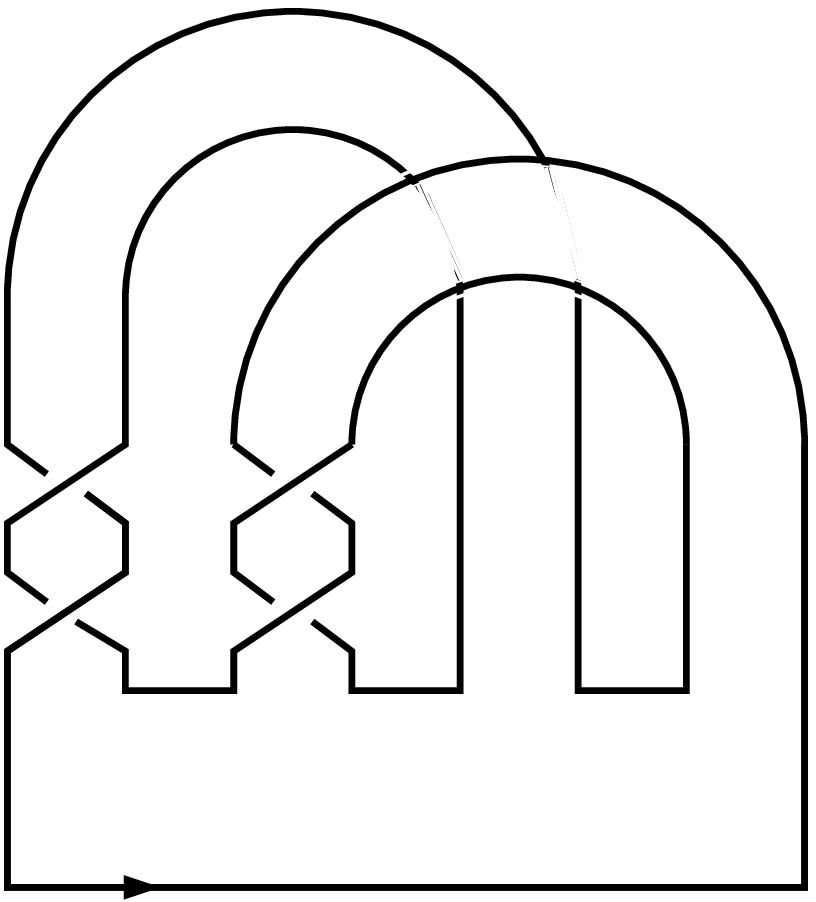,height=38mm,width=38mm,angle=0}
\hspace*{2cm} \psfig{figure=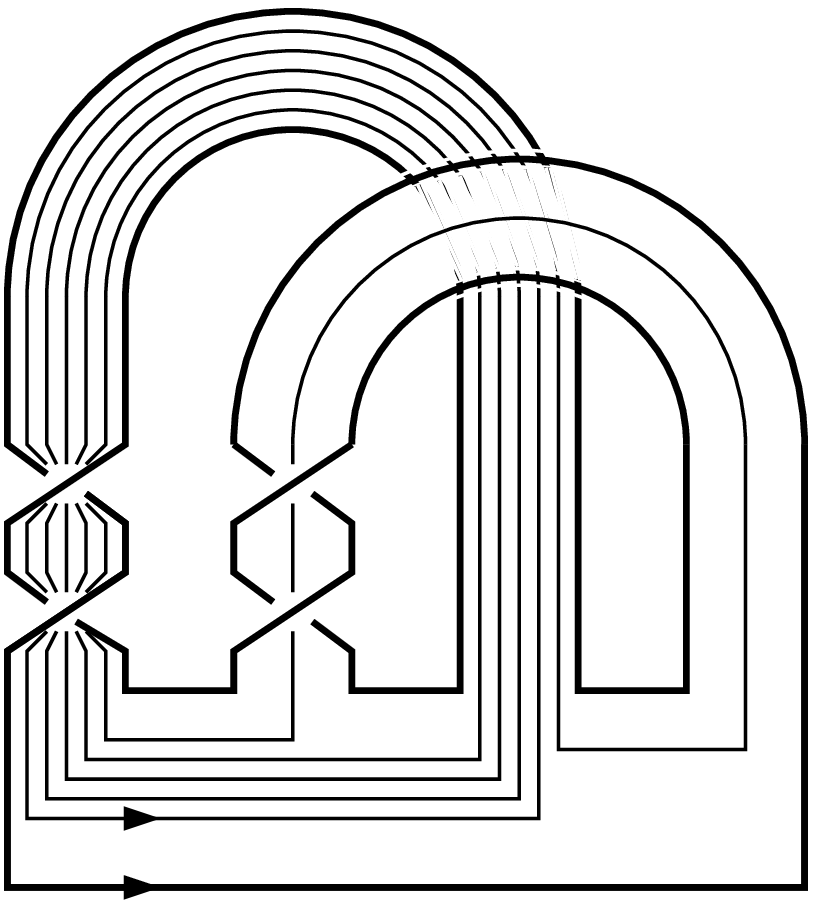,height=38mm,width=38mm,angle=0}} 
\caption{\label{trefoil} {\sl The fiber $\Sigma_{K}$ with the
left-handed trefoil $K$ as boundary (left) and with the knot 
$\gamma_{p} \subset \Sigma_{K}$ (right, with $p=5$). }}
\end{figure}

We point out that the left-handed trefoil knot $K$ is a genus $1$ knot,
whose minimal genus Seifert surface (the fiber $\Sigma_{K}$, 
canonically defined up to isotopy) is a surface having
boundary $K$ itself, illustrated in Figure \ref{trefoil}
(left hand side). The fiber $\Sigma_{L}$ of $L$ is obtained from $\Sigma_{K}$
by removing any disk.

We will now identify now a family of simple closed curves,
lying in $\Sigma_{L}$, that are homologous as elements of
$H_{1}(S^{3} \setminus \nu L,\Z)$, but not isotopic.
We have several choices of simple closed curves lying in $\Sigma_{K}$.
For sake of definiteness, we will
consider the infinite family of knots defined as in the right hand side of 
Figure \ref{trefoil}: the $p$-th element of the family has $p$ strands on the
first twisted annulus and $1$ on the second, as represented in Figure
\ref{trefoil} for the case of $p=5$ (in what follows we will always consider
$p \geq 2$). Note that each $\gp$ has a natural framing induced by the surface
$\Sigma_{K}$.

\begin{figure}[h] 
\centerline{\psfig{figure=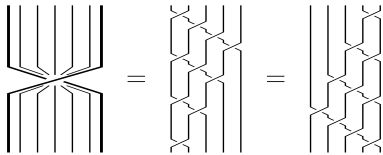,height=20mm,width=50mm,angle=0}}
\caption{\label{braiding} {\sl The effect of the half-twist of the band
on the strands, in two equivalent form.}}
 \end{figure}

\begin{figure}[h] 
\centerline{\psfig{figure=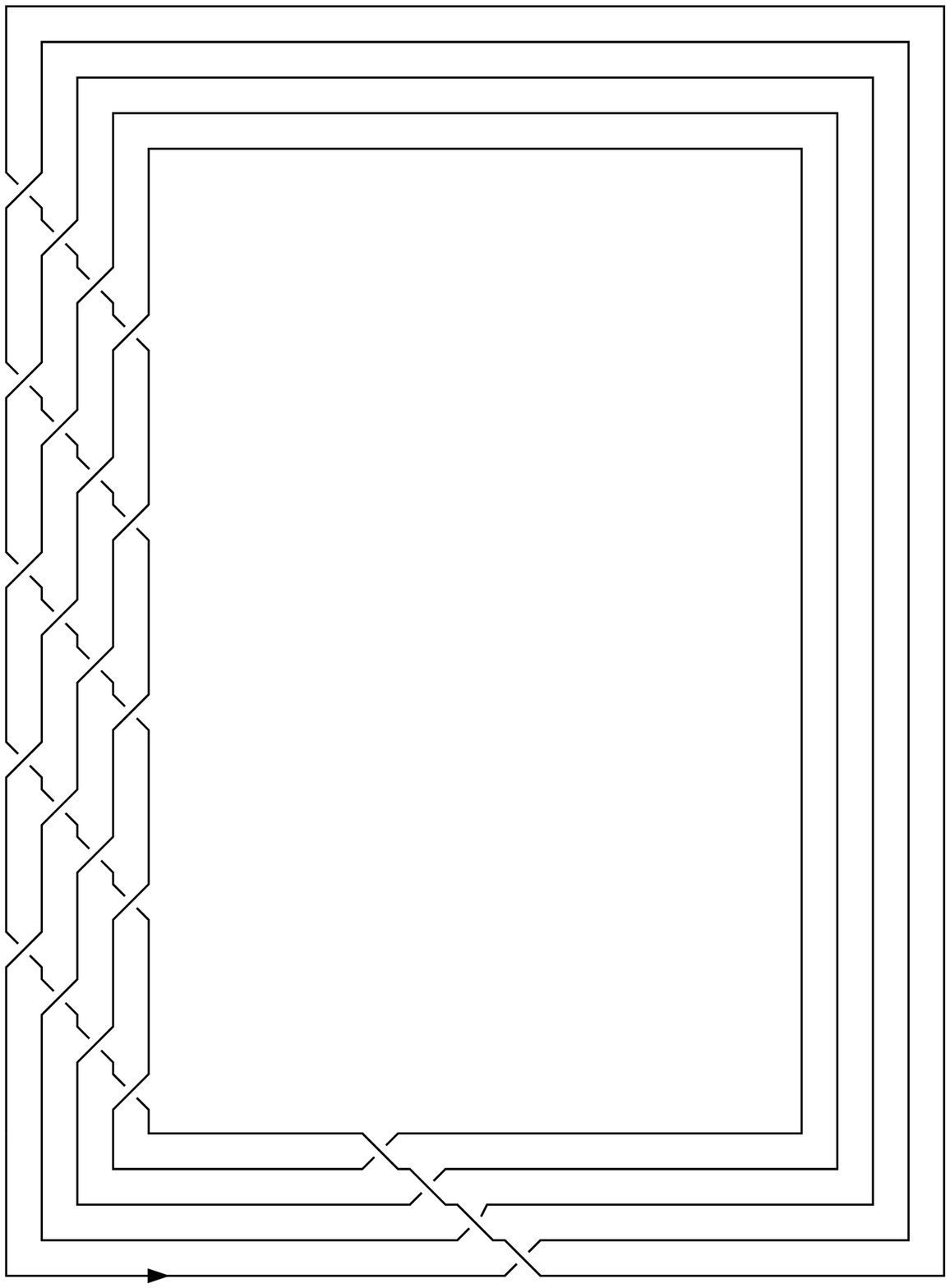,height=44mm,width=33mm,angle=0}
\hspace*{3mm} \psfig{figure=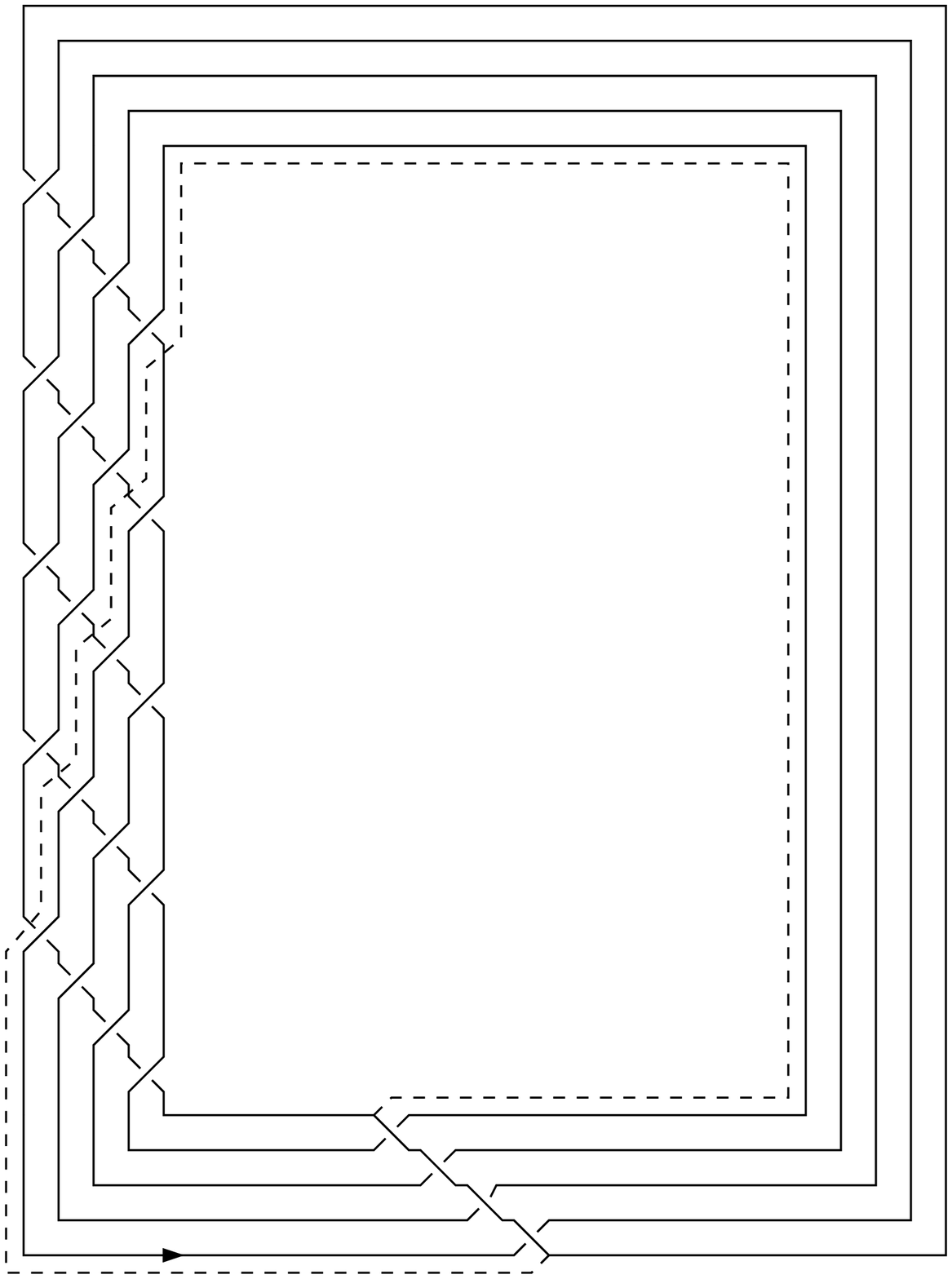,height=44mm,width=33mm,angle=0}
\hspace*{3mm} \psfig{figure=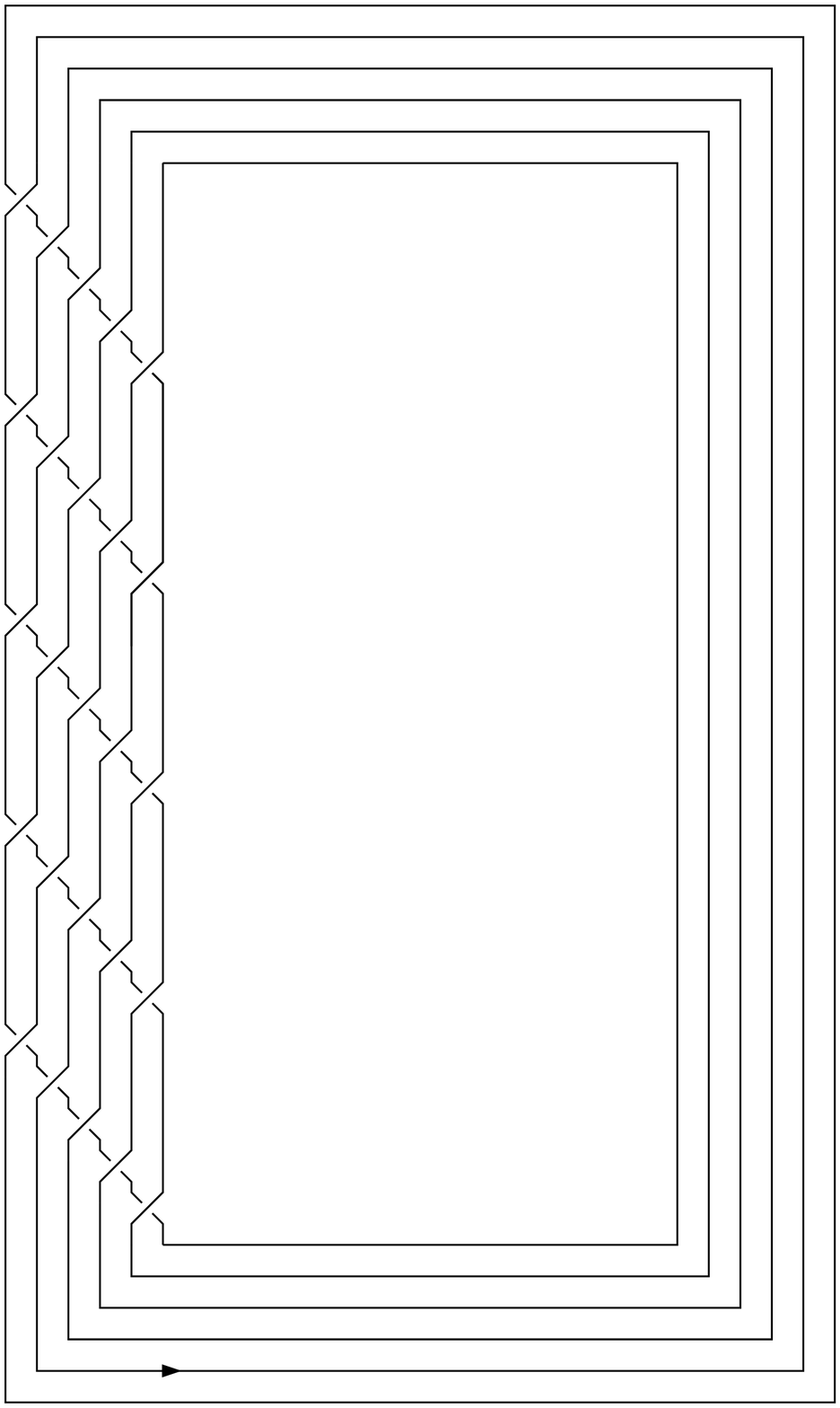,height=44mm,width=33mm,angle=0}}
\caption{\label{steps} {\sl A standard presentation of the torus knot through
the isotopy indicated by the dashed line . }}
\end{figure}

It is not difficult to recognize that the knot denoted by $\gp$ is in 
fact the torus knot $T(p,p+1)$. In fact, the two half-twists on the
annulus act on the $p$ strands as illustrated in Figure \ref{braiding}
(for $p=5$).
Every strand under-crosses all the strands at its right, and we can represent
this with both the braids of Figure \ref{braiding}
(the equivalence of the two braids can be checked with repeated application
of the braid relations). We can now represent $\gp$
as the composition of the two braids of Figure \ref{braiding}, plus
the additional crossing of the $p$-th strand over the others (due to the
passage on the second annulus of $\Sigma_{K}$).
We obtain in this way the presentation of Figure \ref{steps} (left).
Proceeding now as in Figure \ref{steps}, we obtain a presentation of $\gp$
as closure of a braid with $(p+1)$ strands
that shows clearly that it is a torus knot, homologous to $p$ times the meridian
and $(p+1)$ times the longitude of a standard torus in $S^{3}$.
In particular, the example on the right hand side of Figure 
\ref{trefoil} is the torus knot $T(5,6)$.

\begin{figure}[h] 
\centerline{\psfig{figure=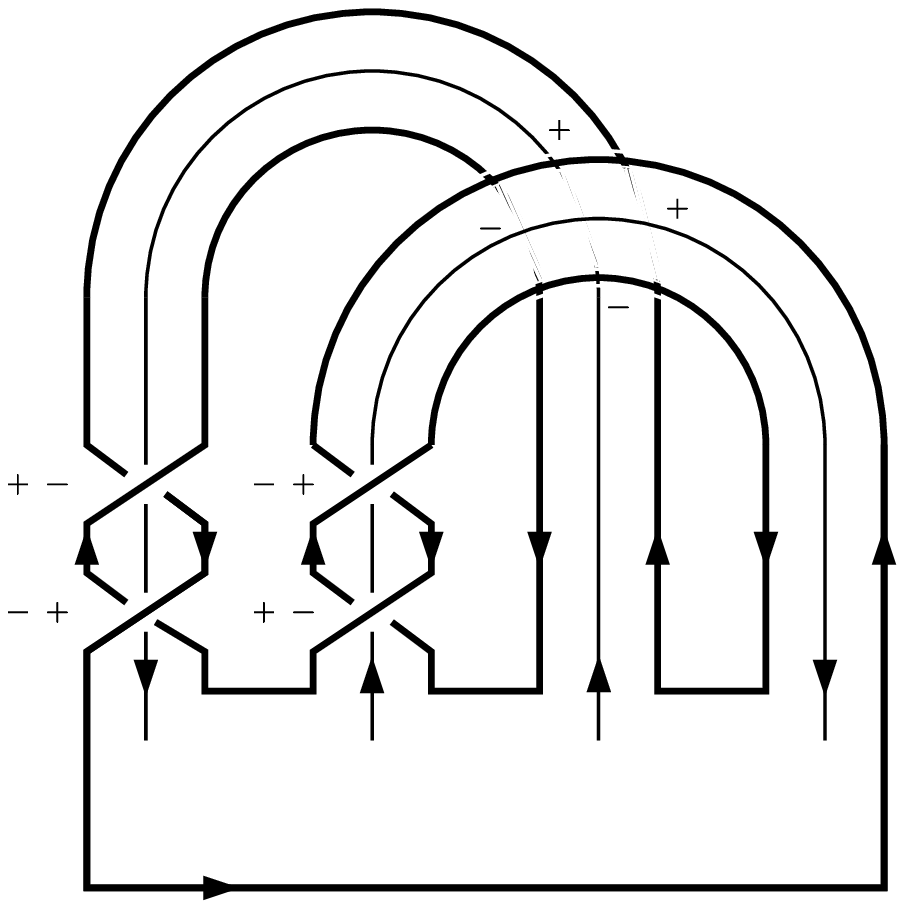,height=38mm,width=38mm,angle=0}
\hspace*{2cm} \psfig{figure=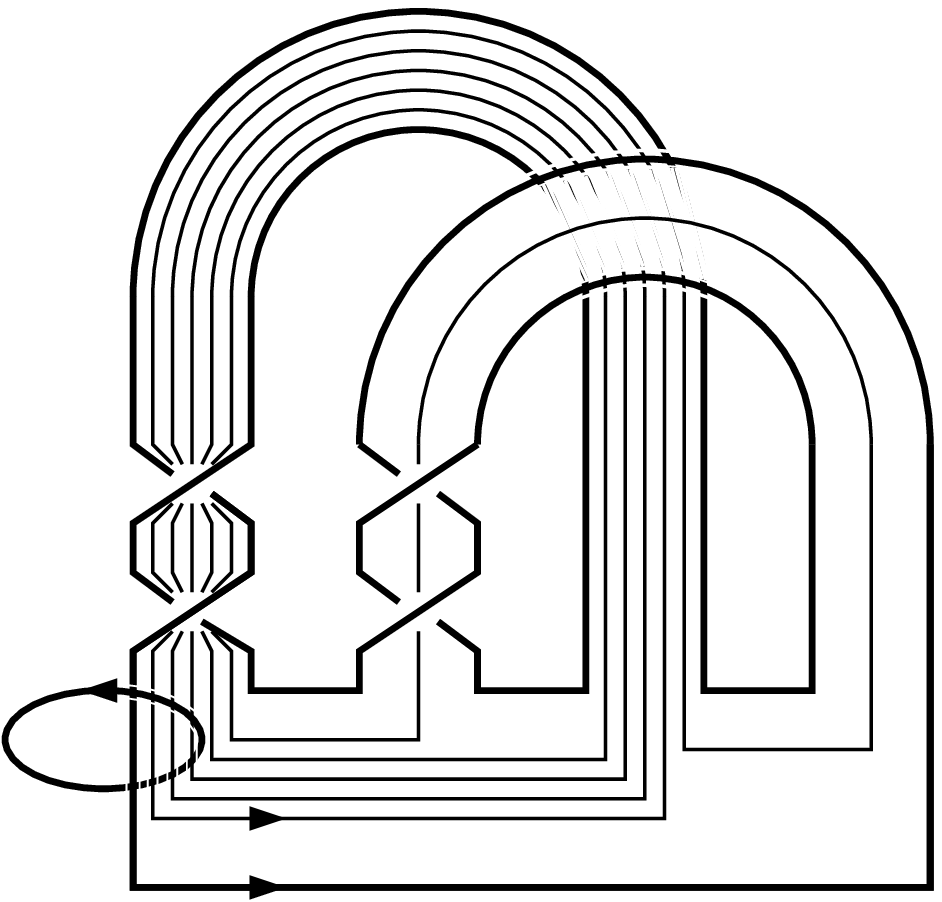,height=38mm,width=44mm,angle=0}}
\caption{\label{crossing} {\sl The crossing numbers of $\gp$ (left) and the 
curve $\gamma_{5,3} \subset \Sigma_{L}$, together with the pair $K \cup M$ 
 (right). }}
\end{figure}

We claim that the linking number of $\gp$ with
$K$ is zero. In fact, by looking at the crossings of $\gp$ and $K$ in
Figure \ref{trefoil}, we can see that
the computation can be reduced to verifying that the sum of the signs
of the crossings of each of the two strands of 
Figure \ref{crossing} (left hand side)
and $K$ equals zero, something that can be checked by
direct computation as indicated in the Figure.

Now let's consider the link $L$ and its fiber $\Sigma_{L} \subset \Sigma_{K}$.
Without loss of generality, we can assume that the curves $\gp$ are
embedded in $\Sigma_{L}$, but in order to do so we need to specify how
$M$ links $\gp$. We will fix an integer $q \geq 1$, and consider the family
of curves $\gamma_{p}$ for $p \geq q$. 
We can isotope $\gp$ (in $\Sigma_{K}$) 
in such a way that $\gp$ pierces with the first $q$ of its strands the 
spanning surface of $M$, as illustrated in the
right hand side of Figure \ref{crossing} (where $p=5$, $q=3$).
We will denote by $\gpd$ the curve of $\Sigma_{L}$ thus obtained.
In particular, for $q = 1$, $M$ is
isotopic (in $S^{3} \setminus \nu \gamma_{p,1}$) 
to a meridian of $\gamma_{p,1}$. As knots in $S^{3}$, we have 
$lk(\gpd,M) = q$.

We want to make clear that, with the definition above, the surface $\Sigma_{L}$
contains, for all values of $q$, the entire 
family of curves $\{\gpd\}_{q,p \geq q}$.

We can define now a $3$-component link 
$L_{p,q} = L \cup \gpd$, which has linking matrix
\begin{equation} \label{lima} l_{p,q} = 
\left( \begin{array}{ccc} - & 1 & 0 \\ 1 & - & q \\ 0 & q & - 
\end{array} \right). \end{equation}
Observe that the linking matrix does not depend on $p$.
As a consequence of this we have, in the homology of 
$S^{3} \setminus \nu L$, the following relation:
\begin{equation} \label{homca}
[\gpd] = lk(\gpd,K) [\mu(K)] + lk(\gpd,M) [\mu(M)] = q [\mu(M)] 
\in H_{1}(S^{3} \setminus \nu L,\Z). \end{equation}  
Although all the $\gpd$'s (for a fixed value of $q$) are homologous,
they are quite clearly non isotopic, for different values of $p$, as knots in
$S^{3}$ and, a fortiori, in $S^{3} \setminus \nu L$. We will exploit this fact in
order to prove Theorem \ref{infinite}.

\section{Construction of the symplectic link surgery manifold} \label{manifold}
The link surgery construction is a convenient method to translate some properties
of knots and links to $4$-dimensional manifolds. In our case
we will construct a symplectic manifold, homotopy equivalent to the elliptic surface
$E(2)$, starting from the link $L$.
Although not strictly necessary for the proof of Theorem \ref{infinite},
we will present the manifold in two different ways. First, as knot surgery
for the knot $K$, and next as link surgery manifold for the link $L$.
This pretty straightforward observation, true for any knot $K$, has some 
interest per se.

We will start with the standard definitions contained in \cite{FS1} of a 
homotopy $E(n)$ associated with the knot $K$: 
\begin{equation} \label{linksu}
E(n)_{K} = (E(n) \setminus \nu F) \cup_{T^{3}} S^{1} \times
(S^{3} \setminus \nu K) = E(n) \#_{F = S^{1} \times C_{K}} S^{1} \times N_{K}
\end{equation}
where the gluing map on the boundary $3$-torus identifies the elliptic fiber
$F$ with $S^{1} \times \mu(K)$ and, reversing 
orientation, $\lambda(K)$ with the meridian circle to the elliptic fiber.
In Equation \ref{linksu} the $3$-manifold $N_{K}$ is the result of
$0$-surgery of $S^{3}$ along $K$ and $C_{K}$ is the core
of the Dehn filling (a curve isotopic to $\mu(K)$). The torus 
$S^{1} \times C_{K}$ has a canonical framing induced by the Dehn filling in
$N_{K}$, while $F$ in $E(n)$ has a canonical framing induced from the
fibration of $E(n)$. The fiber sum identifies the two tori and acts as 
complex conjugation on the normal bundles.
We are interested in the case of $n=2$:
we claim that $E(2)_{K}$ can be described as link surgery manifold
over the link $L$, with the generalized definition introduced in \cite{V1}.
Precisely, we can construct an homotopy $E(2)$ starting from $L$ and the
fibration of $S^{3} \setminus \nu L$ of fiber $\Sigma_{L}$ in the following
way: \begin{equation} \begin{array}{c} \label{othdef} E(2)_{K} =
\coprod_{i=1}^{2} (E(1)_{i} \setminus
\nu F_{i} ) \cup_{T^{3}_{i}} S^{1} \times (S^{3} \setminus \nu L) = \\ \\
= E(1)_{1} \#_{F_{1} = S^{1} \times C_{K}} S^{1} \times N_{K} \#_{F_{2} = 
S^{1} \times C_{M}} E(1)_{2}. \end{array} \end{equation}
The gluing maps identify, on the first $3$-torus, 
$F_{1}$ with $S^{1} \times \mu(K)$ and the meridian of $F_{1}$ with 
$- \lambda(K)$ and, on the second $3$-torus, $F_{2}$ with $S^{1} \times
\lambda(M)$
and the meridian of $F_{2}$ with $\mu(M)$. $C_{K}$ is again the core of the
Dehn filling along $K$ and $C_{M}$ is the core of the trivial Dehn filling
of the $\infty$-surgery along $M$ (a curve isotopic to $\lambda(M)$ itself).
Note that, with this gluing prescription above, 
$\Sigma_{L}$ is naturally capped off with one disk section in
each $(E(1) \setminus \nu F_{i})$. As $\lambda(M)$ is isotopic to $\mu(K)$
we have, up to isotopy, $C_{K} = C_{M}$ and  
$F_{1} = F_{2} = S^{1} \times C_{(\cdot)}$.
The connection between the two definitions above comes by observing that, as 
$E(2) = E(1) \#_{F} E(1)$, we can think to $E(2)_{K}$ as obtained by
fiber summing $E(1)$ to $E(1)_{K}$ along a copy of $F = S^{1} \times C_{K}$, 
which is exactly the definition of Equation \ref{othdef}. Note that the
definitions of Equations \ref{linksu} and \ref{othdef} are equivalent also
in the symplectic category. 

The definition of $E(2)_{K}$ in Equation \ref{othdef} shows the existence of
two noteworthy homology classes of $E(2)_{K}$, images under the injective map
\begin{equation} i_{*}: H_{1}(S^{3} \setminus \nu L,\Z) \longrightarrow
H_{2}(S^{1} \times (S^{3} \setminus \nu L),\Z) \longrightarrow
H_{2}(E(2)_{K},\Z) \end{equation} of the two generators $[\mu(K)]$ and
$[\mu(M)]$.
The first one is identified with the homology class of $F$, while the second 
is the homology class of a kind of rim torus,
identified with $S^{1} \times \mu(M)$, that we will denote by $R$. This class 
is primitive: in order to verify this, observe that the torus $R$ intersects
once (with positive sign) the sphere obtained by capping the annulus spanning
$M$ and pierced once by $K$ (representing the dual to the class 
$(0,1) \in H^{1}(S^{3} \setminus \nu L,\Z)$) with one vanishing disks in each
copy of $E(1) \setminus \nu F$. The class of $R$ has self-intersection $0$.

\section{Infinitely many lagrangian tori}
We are ready now, using the constructions of Sections \ref{linkco} and 
\ref{manifold}, to exhibit a family of lagrangian tori that represents
any multiple of the class of the rim torus $R$. Fix, as in Section 
\ref{linkco},
an integer $q \geq 1$ and pick the collection of curves $\gpd$
with $p \geq q$. The following holds true.
\begin{lm} \label{homol} The images $R_{p,q} \subset E(2)_{K}$ of the tori
$S^{1} \times \gpd \in S^{1} \times (S^{3} \setminus \nu L)$
define a family of lagrangian, framed tori 
representing the class $q[R]$. \end{lm}
{\bf Proof:} The homology class of these tori is given, by Equation 
\ref{homca}, 
by the formula \begin{equation} \begin{array}{c} [R_{p,q}] = i_{*} [\gpd] = lk(\gpd,K)
i_{*} [\mu(K)] + lk(\gpd,M) i_{*}[\mu(M)] = \\ \\ = q[R] \in H_{2}(E(2)_{K},\Z).
\end{array} \end{equation}
The elements of the
infinite collection of tori $R_{p,q}$, with $p \geq q$,
is therefore homologous to $q$ copies of $R$.
As the knot $K$ is fibered, the closed $3$-manifold
$N_{K}$ is fibered too, with fiber the surface $\Sigma_{L}$ capped off with 
two disks. The fiber sum presentation of Equation \ref{othdef},
scaling suitably the symplectic forms to obtain matching volumes on the
gluing tori, shows that $E(2)_{K}$ admits a symplectic structure that, 
on the interior of $S^{1} \times N_{K} \setminus \nu (S^{1} \times C_{K} \cup 
S^{1} \times C_{M}) = S^{1} \times (S^{3} \setminus \nu L)$, coincides with
the restriction of a standard symplectic
form on $S^{1} \times N_{K}$, namely $\alpha \wedge dt + \epsilon \beta$,
where $\alpha$ is a closed $1$-form on $N_{K}$ determined up to isotopy by
the fibration of $N_{K}$ and $\beta$ is a closed $2$-form on $N_{K}$
restricting to a volume form on each fiber. The tori $S^{1} \times \gpd$
embed in $S^{1} \times N_{K}$ and are lagrangian with respect to its
symplectic structure, as the tangent space in each point
is spanned by $\frac{\partial}{\partial t}$ and by a vector tangent to the
fiber, and therefore lying in the kernel of $\alpha$. Consequently, the tori
$R_{p,q}$ are lagrangian and have a canonical framing induced from
$\Sigma_{L}$; this framing coincides with the lagrangian framing, as the tori
obtained by pushing off $\gamma_{p,q}$ along $\Sigma_{L}$ are lagrangian.
This completes the proof of the Lemma. {\it Q.e.d.} \endproof

Lemma \ref{homol} asserts that the tori $R_{p,q}$, for a fixed value of $q$,
are images of homotopic lagrangian embeddings of $T^{2}$ in $E(2)_{K}$.
Our goal is now to show that these tori are not isotopic. In order to do so,
we can define the family of (symplectic) manifolds
\begin{equation} X_{p,q} = E(1) \#_{F=R_{p,q}} E(2)_{K}. \end{equation}
The definition of the fiber sum depends a priori on the choice of the framing
for $R_{p,q}$, determining (with a marking of the tori) the gluing map for the
boundary $3$-tori. Anyhow, any orientation preserving 
self-diffeomorphism of $\partial(E(1) \setminus \nu F)$ extends to $(E(1)
\setminus \nu F)$ (see \cite{GS}): different choices of the framing or 
marking lead therefore to the same smooth manifold, namely $X_{p,q}$
depends only on the isotopy class of $R_{p,q}$. More is actually true,
namely the smooth structure of $X_{p,q}$ depends only the diffeomorphism type of
the pair $(E(2)_{K},R_{p,q})$.
Our goal is to distinguish different smooth structures,
for different values of $p$, by computing the Seiberg-Witten polynomial of
$X_{p,q}$. In order to do so, we start with the following
\begin{pr} \label{samman} Consider the manifold
$ X_{p,q} = E(1) \#_{F=R_{p,q}} E(2)_{K}$ and the
$3$-component link $L_{p,q} = K \cup M \cup \gpd$.
The manifold $X_{p,q}$ is diffeomorphic to the link surgery manifold
$E(3)_{L_{p,q}}$. \end{pr}
{\bf Proof:} The manifolds $X_{p,q}$ and $E(3)_{L_{p,q}}$ can both be written,
as smooth manifolds, as
\begin{equation} \coprod_{i=1}^{3} (E(1) \setminus \nu F) \cup_{T^{3}_{i}}
S^{1} \times (S^{3} \setminus L_{p,q}), \end{equation} where the gluing maps on the boundary
$3$-tori are defined in different ways for the two manifolds.
However, because of the aforementioned extension property of the
self-diffeomorphisms of
$\partial(E(1) \setminus \nu F)$, the choice of the gluing maps does not
affect the smooth structure of the resulting manifold. {\it Q.e.d.} \endproof
Proposition \ref{samman} allows us to apply the following Lemma of
Fintushel-Stern (\cite{FS1}).
\begin{lm} \label{ronron}
Let $\Delta^{s}_{p,q}(x,y,t)$ be the symmetrized Alexander
polynomial of
the link $L_{p,q} = K \cup M \cup \gpd$. Then the Seiberg-Witten
polynomial of $X_{p,q}$ is given by
\begin{equation} \label{relsw} SW_{X_{p,q}} = \Delta^{s}_{p,q}(x^{2},y^{2},t^{2})
\end{equation}
where we identify $x,y,t$ with (the Poincar\'e dual of) the images of
$S^{1} \times \mu(K), S^{1} \times \mu(M),
S^{1} \times \mu(\gpd)$ respectively. \end{lm}
Note that, as the relative Seiberg-Witten polynomial of $(E(1),F)$ is
just equal to $1$, the polynomial of Equation \ref{relsw} is in fact the
relative Seiberg-Witten polynomial of $(E(2)_{K},R_{p,q})$ - see \cite{McMT} for
a discussion of the invariants from this point of view.
We will use this result in the next section to prove our main Theorem.
\section{Proof of the main theorem}
In order to prove Theorem \ref{infinite} we could attempt to compute,
using Lemma \ref{ronron}, the SW polynomial of $X_{p,q}$ for each value of $p$
and use the result
to distinguish the smooth structure for different values of $p$
and, with that, the isotopy class of the torus $R_{p,q}$.
Conceptually, the computation of the Alexander polynomial of $L_{p,q}$
does not present any difficulty, but obtaining a general formula
is practically not viable. Moreover, even when this is
done, while comparing the SW polynomial of two manifolds, we would need to
show that there is no change of basis in the second cohomology group
transforming one polynomial in the other, compare with the discussion in
\cite{V}.

For this reason, we will be content with a weaker result, that is anyhow
sufficient to prove Theorem \ref{infinite} (even in its stronger form for
diffeomorphisms of pairs)
\begin{thr} \label{lot} For any choice of $q \geq 1$
the family of manifolds $\{X_{p,q}\}_{p \geq q}$ contains infinitely many
pairwise non-diffeomorphic manifolds. \end{thr}
{\bf Proof:} In order to prove the statement it is enough to show
that the number of basic classes of $X_{p,q}$, that we will denote by 
$\beta_{p,q}$, satisfies the condition $lim_{p}\beta_{p,q} = \infty$. 
As all the manifolds $X_{p,q}$ have a finite number of basic classes, 
this implies that infinitely
many components of the family $\{X_{p,q}\}_{p \geq q}$ are distinguished by the SW 
polynomial. To start we can observe that $\beta_{p,q}$ 
is the same as the number of nonzero terms of
$\Delta_{p,q}(x,y,t)$
(symmetrization is irrelevant here).
The latter number, moreover, is bounded below by the number of nonzero
terms of any specialization of the polynomial, so we will make our life 
simpler by considering specializations of $\Delta_{p,q}(x,y,t)$. We have the
following Lemma.
\begin{lm} The number of nonzero terms of the Alexander polynomial 
$\Delta_{p,q}(x,y,t)$ of 
$L_{p,q}$ is bounded below by the number of nonzero terms of the
polynomial (written here in quotient form)
\begin{equation} \label{nepol} P_{p,q}(t) = 
\frac{(1-t^{q})(1-(t^{p+1})^{p})}{(1-t^{p})(1-t^{p+1})}. \end{equation} \end{lm}
{\bf Proof:} The number of nonzero terms of $\Delta_{p,q}(x,y,t)$
is bounded below by the number of nonzero terms of $\Delta_{p,q}(1,y,t)$.
This specialization of the Alexander polynomial can
be computed using the Torres formula (see e.g., \cite{Tu})
to get \begin{equation} \label{special} \Delta_{p,q}(1,y,t) =
(y^{lk(M,K)}t^{lk(\gpd,K)}-1) \Delta_{\gpd \cup M}(y,t) =
(y-1) \Delta_{\gpd \cup M}(y,t). \end{equation}
Let's consider first the case of $q=1$. In this case, $\gamma_{p,1} \cup M$
is a pair given by the torus knot and its meridian or, which is the same, its
$(0,1)$-cable. As such we
can represent it as result of splicing $\gp$ with the generalized Hopf link 
(the ``necklace'') $H$ (with $3$ components), along a component $H_{2}$. 
The Alexander polynomial of $H$ is given, as easily checked, 
by $\Delta_{H}(u,v,w) = (u-1)$. Applying the splicing formula of Theorem 5.3 
of \cite{EN} we find that $\Delta_{\gamma_{p,1} \cup M}(y,t) = 
\Delta_{\gp}(t)$, where the latter is the 
Alexander polynomial of the torus knot $T(p,p+1)$, namely
\begin{equation} \label{torpol} \Delta_{\gp}(t) = \frac{(1-t)(1-(t^{p+1})^{p})}{(1-t^{p})(1-t^{p+1})}.
\end{equation}
This exactly the polynomial $P_{p,1}(t)$ of Equation \ref{nepol} and it should
be clear at this point, by looking at Equation \ref{special}, that 
the statement of the Lemma holds true for $q=1$.
To prove the general case, let's write 
\begin{equation} \label{coeff} \Delta_{\gpd \cup M}(y,t) = \sum_{k} 
a_{k,p,q}(y)t^{k},
\end{equation} where the
Laurent polynomials $a_{k,p,q}(y)$ are defined by this identity.
We can write Equation \ref{special} in the form
\begin{equation} \Delta_{p,q}(1,y,t) = \sum_{k} (y-1)a_{k,p,q}(y)t^{k}. 
\end{equation}
Consider the set of coefficients $(y-1)a_{k,p,q}(y)$ of this Laurent polynomial
of the variable $t$. The number of 
coefficients that are nonzero is bounded below by the number of terms
$a_{k,p,q}(1)$ that are nonzero. By the definition contained in Equation 
\ref{coeff}, 
the value of $a_{k,p,q}(1)$ is determined by the Equation
\begin{equation} \sum_{k} a_{k,p,q}(1)t^{k} =  \Delta_{\gpd \cup M}(1,t)
= \frac{t^{lk(\gpd,M)}-1}{t-1}\Delta_{\gp}(t) =  
\frac{t^{q}-1}{t-1}\Delta_{\gp}(t),
\end{equation} where to get the last expression we used again Torres formula.
By looking at Equation \ref{torpol} we see that this last polynomial is 
exactly the polynomial $P_{p,q}(t)$ of Equation \ref{nepol}. {\it Q.e.d.} \endproof

The formula of Equation \ref{nepol} tells us that $\beta_{p,q}$ is bounded
below by the number of nonzero terms in the polynomial $P_{p,q}(t)$.
Concerning this polynomial,
we have the following result:
\begin{lm} \label{est} 
The number of nonzero terms of the polynomial $P_{p,q}(t)$, with $p \geq q$, 
is bounded below by $p-q+1$. \end{lm}
{\bf Proof:} The proof of this statement follows by a fairly simple argument.
First, we can rewrite, in $\Z[[t]]$,
\begin{equation} \label{longpoly}\begin{array}{c}
P_{p,q}(t) = \frac{(1+t^{p+1}+ ... + 
t^{(p-1)(p+1)})(1-t^{q})}{1-t^{p}} = \\ \\ = \frac{1-t^{q}+t^{p+1}-t^{p+1+q}+ ... + t^{(p-1)(p+1)}
- t^{(p-1)(p+1) + q}}{1-t^{p}} =
\\ \\ = (1-t^{q}+t^{p+1}-t^{p+1+q}+ ... + t^{(p-1)(p+1)}
- t^{(p-1)(p+1) + q})(1+t^{p} + t^{2p} + ... ).
\end{array} \end{equation}
The product of the polynomial and the formal power series above is in fact a 
polynomial, and we claim that it contains  
all the $p-q+1$ terms of the form \begin{equation} 
\label{terms} t^{(q-1)(p+1)},t^{(q-1)(p+1)+p},...,t^{(q-1)(p+1)+ (p-q-1)p},
t^{(q-1)(p+1)+ (p-q)p}. 
\end{equation} In order to prove so, observe first that, as $0 \leq (q-1) 
\leq (p-1)$, 
the polynomial in the second line of Equation \ref{longpoly} contains
the term $t^{(q-1)(p+1)}$ with coefficient equal to $1$.
Therefore we obtain, with coefficient equal to $1$, all the terms of 
Equation \ref{terms} by multiplying the term $t^{(q-1)(p+1)}$ by the 
first $(p-q+1)$ terms of the formal power series (plus an infinite series of
other terms that, in fact, will get canceled). 
Next, considering the fact that the power series has all
terms with positive coefficients, the only possibility for the terms 
$t^{(q-1)(p+1) + np}$,
$0 \leq n \leq (p-q)$
to be canceled is that there exists a pair $k,l \in \N$,
with $0 \leq l \leq (p-1)$, such that
\begin{equation} \label{coeffeq} (q-1)(p+1)+np = kp + l(p+1) +q \end{equation}
(roughly speaking, the power of $t^{l(p+1)+q} \cdot t^{kp}$, appearing
with a negative sign, must equal the power of $t^{(q-1)(p+1) + np}$).
Remember that, in this equation, $p$ and $q$ are fixed.
Equation \ref{coeffeq} means that, {\it mod} $p$, we must have
$l + 1 = 0$. As $0 \leq l \leq (p-1)$, the only solution to this condition
is $l =(p-1)$. With this value of $l$, Equation \ref{coeffeq} becomes
\begin{equation} (q-1)p + np = kp + p^{2}. \end{equation}
The smallest value of $n$ for which this equation holds true is when $k=0$,
for which $n=p-q+1$, namely none of the terms in 
Equation \ref{terms} gets canceled, as claimed (while the higher terms get
canceled, as they should). With similar considerations
it is possible to show that the coefficient of these terms is exactly
equal to $1$, but we will omit the proof
(as it has no implications to our result). {\it Q.e.d.}  \endproof
(Note that the estimate contained in Proposition \ref{est} is quite rough,
and probably without much effort a precise value could be obtained, but
this result would be irrelevant in our discussion).

Proposition \ref{est} implies that $\lim_{p \rightarrow \infty} \beta_{p} 
\geq \lim_{p \rightarrow \infty} (p-q+1) = \infty$, and thus completes the 
proof of Theorem \ref{lot} and, with that, Theorem \ref{infinite}. {\it Q.e.d.}

\end{document}